\documentclass{article}[12pt]

\usepackage{graphicx}
\usepackage{tabu}
\usepackage{epsfig}
\usepackage{mathtools}
\usepackage{float}
\usepackage{algorithm}
\usepackage[noend]{algpseudocode}
\usepackage{amsmath}
\usepackage[export]{adjustbox}
\usepackage{subcaption}
\usepackage{tikz}
\usepackage{amsmath}

\newcommand{\expt}[1]{\langle#1\rangle}

\usepackage{amssymb}
\usepackage{authblk}

 \title{Why Simple Quadrature is just as good as Monte Carlo}

 \author{Kevin Vanslette, Abdullatif Al Alsheikh, and Kamal Youcef-Toumi }
 \affil{\small{\emph{ Massachusetts Institute of Technology}}}
  \affil{\small{\emph{King Abdulaziz city for Science and Technology}}}

 \date{\today}

\begin{document}
\maketitle

\abstract{We motive and calculate Newton--Cotes quadrature integration variance and compare it directly with Monte Carlo (MC) integration variance.  We find an equivalence between deterministic quadrature sampling and random MC sampling by noting that MC random sampling is statistically indistinguishable from a method that uses deterministic sampling on a randomly shuffled (permuted) function.  We use this statistical equivalence to regularize the form of permissible Bayesian quadrature integration priors such that they are guaranteed to be objectively comparable with MC. This leads to the proof that simple quadrature methods have expected variances that are less than or equal to their corresponding theoretical MC integration variances. Separately, using Bayesian probability theory, we find that the theoretical standard deviations of the unbiased errors of simple Newton--Cotes composite quadrature integrations improve over their worst case errors by an extra dimension independent factor $\propto N^{-1/2}$. This dimension independent factor is validated in our simulations.}

\numberwithin{equation}{section}
\newtheorem{thm}{Condition}[section]

\section{Introduction}
    %
Some of the most predominant computational techniques for estimating intractable integrals numerically are based on Monte Carlo (MC) random sampling \cite{Metro,GA,SS,Spectral,Dakota,stan}. These methods aim to combat the infamous ``curse of dimensionality" and allow us to reliably access a wider range of computationally and theoretically challenging problems.  MC integration has thus positively influenced fields of physics, engineering, computer science, statistics, and machine learning.

The fidelity of MC integration and Newton--Cotes quadrature integration are only ever loosely compared because MC is random and uses standard deviations while quadrature is deterministic and uses worst case errors. One way to gauge the fidelity of an integration estimate is by how its standard deviation or error improves with the number, $N$, of $D$ dimensional samples of the integral. For simple MC, the standard deviation of the integral estimate improves $\propto N^{-\frac{1}{2}}$, which is independent of $D$. Letting $n_{\mathcal{Q}}$ be the number of derivatives taken in the Newton-Cotes quadrature error method $\mathcal{Q}$, quadrature error has poor performance, $\propto N^{-\frac{n_{\mathcal{Q}}}{D}}$, when $D$ is large.  Due to this dependence on the dimension, an indirect comparison between MC and quadrature suggests that MC estimation is better suited for high dimensional integrals, which is used to motivate MC integration methods throughout the majority of the literature. This article offers a different view. 

 We believe the main difficulty underpinning the lack of direct comparisons between the fidelity of quadrature and MC actually originates from the lack of reconciliation between Bayesian and frequentist probability theory. Bayesian probability distributions represent a \emph{user's} informational state about a given uncertain system. The probabilistic description of the system is therefore subject to the user's prior knowledge and it is designed to change when new information is learned to best reflect the user's new informational state \cite{Catichaupdate,MyThesis}. Quadrature integration variances are typically Bayesian in nature \cite{MCOHagan,Briol}. In contrast, theoretical frequentist probability distributions and expectation values are treated as fixed values associated to the physical random system rather than being associated to its user's knowledge of said system. These values may be realized with certainty from the population statistics in the infinite sample limit, which is the basis of MC. Thus, for a given informational state, Bayesian and frequentist inferences may actually disagree with one another -- even to the extent of not agreeing on the set of possible outcomes -- if they are not properly regulated with respect to one another. Thus, if not regulated, it is difficult to construct fair comparisons between these methods. 
 

In this article we motive and calculate simple quadrature variances and compare them directly to MC by proposing and implementing a Bayesian to frequentist regularization called \emph{counterfactual shuffling} (CS).  First, we find that MC sampling an integrand function is equivalent in probability to simple quadrature sampling a randomly shuffled integrand function. This allows us to regularize the Bayesian probability distributions associated with quadrature by imposing that \emph{if} the integrand function were shuffled, then its resulting Bayesian probability distribution \emph{should} match its frequentist MC counterpart -- which is the CS constraint. Using the law of total variance and under the CS constraint, when there is prior information available about the integrand function, we find that the expected variance of a simple quadrature procedure is less than or equal to the variance of MC.  We make comments about CS in Bayesian and frequentist methods, the difficulty of merging frequentist and Bayesian frameworks, and on error estimation in estimation theory.  



Secondly, we improve the Newton--Cotes composite quadrature error analysis using Bayesian probability. 
We recognize that in practice the error of any estimate is uncertain. This allows us to treat the sum of composite quadrature errors as a set of Bayesian random variables and find that the standard deviation of the unbiased error is theoretically $\propto N^{-(\frac{n_{\mathcal{Q}}}{D}+\frac{1}{2})}$ for the rectangular and midpoint rule. The standard deviation converges faster than the error by an extra factor of $N^{-\frac{1}{2}}$. We simulate and confirm that this extra $N^{-\frac{1}{2}}$ persists in high dimensions. 


The remainder of the article is organized as follows. In Section \ref{background} we review the fidelity of some simple numerical integration methods. In Section \ref{quadfid} we formulate CS and use it to make direct comparisons between quadrature and MC estimation variances. In Appendix \ref{appendix2} we briefly compare quadrature to other MC based methods like MCMC, importance sampling, stratified sampling, and Latin hypercube when possible. In Section \ref{sectionerror} we calculate and simulate the unbiased error variances for quadrature error, which improves the error analysis given by the Newton--Cotes formulas. Our simulations are discussed and tabulated in Appendix \ref{appendix}.

\section{Review: Simple Methods of Numerical Integration\label{background}}

Although a review of numerical integration methods can be found in a number of texts, we review it here for notational consistency while emphasizing certain elements for later reference.

Numerical integration methods aim to perform fast and reliable estimates to integrals of the form,
\begin{align}
I=\int_{\vec{x}\in\Omega}f(\vec{x})\,d\vec{x}.\label{int}
\end{align}
The integral is over the $D$ dimensional volume $\Omega$ that contains vectors of the form $\vec{x}=(x_1,...,x_D)$, which has real valued components, and measure $d\vec{x}=\prod_{i=1}^Ddx_i$. The integral can be transformed such that the volume $\Omega=1$ is unitless and spans the unit $D$ cube $[0,1]^{D}$, without loss of generality. It is assumed that the integrals are well behaved, i.e., $|f(\vec{x})|<\infty$.

MC interprets the integral (\ref{int}) as an expectation value over the unit uniform distribution $\rho(\vec{x})=\frac{1}{\Omega}=1$,
\begin{align}
I=\int_{\vec{x}\in\Omega}\,\frac{1}{\Omega}f(\vec{x})\,d\vec{x}=\expt{f},\label{int2}
\end{align}
such that one can make use of statistical methods. By generating $N$ independent and identically distributed (i.i.d.) random samples $\{\vec{x}_i\}=\{\vec{x}_1,..,\vec{x}_N\}_{MC}$ from $\rho(\vec{x})$, evaluating $f$ at these samples to obtain $\{f_i\}=\{f_1,...,f_N\}_{MC}$, and averaging the results, MC approximates (\ref{int}). That is, the MC estimation of (\ref{int}) is given by,
\begin{align}
\hat{I}_{\mbox{\tiny{MC}}}=\frac{1}{N}\sum_{i=1}^Nf(\vec{x}_i)=\frac{1}{N}\sum_{i=1}^Nf_i.\nonumber
\end{align}
Because the samples in $\vec{x}$ are i.i.d. random numbers, so are its evaluations $\{f_i\}_{MC}$, and therefore the standard deviation for each $f_i$ is homogeneous $\sigma_{f_i}=\sigma_{f}$. The theoretical standard deviation of the MC estimate is therefore,
\begin{align}
\sigma_{\hat{I}_{\mbox{\tiny{MC}}}}=\sqrt{\expt{(\frac{1}{N}\sum_{i=1}^Nf_i)^2}-\expt{\frac{1}{N}\sum_{i=1}^Nf_i}^2}=\frac{1}{N}\sqrt{\sum_{i=1}^N\sigma^2_{f_i}}=\frac{\sigma_{f}}{\sqrt{N}},\label{stdMC}
\end{align}
which decreases $\propto N^{-\frac{1}{2}}$ and without an exponential dependence on the integral's dimension $D$. As $N$ goes to infinity, the estimation becomes Gaussian distributed about $I$ with standard deviation (\ref{stdMC}) due to the central limit theorem. This estimate is unbiased. 

As we are interested in comparing quadrature to MC in the context of Bayesian and frequentist probability theory, we place emphasis on the standard deviation $\sigma_f$ used in MC, which is,
\begin{align}
\sigma_{f}^2=\expt{f^2}-\expt{f}^2=\int_{\vec{x}\in\Omega}\,\frac{1}{\Omega}f(\vec{x})^2\,d\vec{x}-\Big(\int_{\vec{x}\in\Omega}\,\frac{1}{\Omega}f(\vec{x})\,d\vec{x}\Big)^2.\label{2}
\end{align}
There is a known practical condition on equation (\ref{2}) that we will denote as condition \ref{def}:
\setcounter{thm}{3}
\begin{thm}
The standard deviation $\sigma_f$ is a definite yet unknown quantity in practice because it is in terms of the fixed integral we are trying to estimate and over a function that we may not know with complete certainty.\label{def} 
\end{thm}
Due to \ref{def}, $\sigma_{\hat{I}_{\mbox{\tiny{MC}}}}$ is a definite yet unknown quantity as well, which implies the MC formalism can only offer weak theoretical measures of fidelity of an individual MC estimate before sampling as is typical in frequentist methods. Nevertheless, expressing the fidelity of MC in terms of $\sigma_f$ gives useful information about its theoretical dependence on $N$ in (\ref{stdMC}).


Composite quadrature methods estimate integrals by evaluating the integrand function at $N$ definite positions and estimating the volume under (or near) the evaluated points in $\Omega$. In this article we will only consider uniform $D$ dimensional sampling grids. Evaluating $f(\vec{x})$ at the grid locations generate a set of evaluations $\{f_1,...,f_N\}$ that populate the $D$ dimensional sampling grid. Quadrature based integration methods may be decomposed into the product of a generalized height $\mathcal{Q}(i,\{f_1,...,f_N\})$ times a generalized width $\Omega_{i}$, which is summed over the $N'\leq N$ widths,
\begin{align}
\hat{I}_{\mathcal{Q}}= \sum_{i=1}^{N'}\mathcal{Q}(i,\{f_1,...,f_N\})\cdot \Omega_{i}.
\end{align}
The widths $\Omega_{i}\in \Omega$ are $D$ dimensional volumes that partition $\Omega$. Further, the estimate can be written as a linear combination of function evaluations,
\begin{align}
\hat{I}_{\mathcal{Q}}= \sum_{i=1}^{N}q_if_i,\label{linearquad}
\end{align}
which is the Newton--Cotes formula, where $\{q_i\}$ are the weights and $\sum_{i=1}^{N}q_i=N'$. A nonexhaustive set of quadrature methods are: rectangular, midpoint, trapezoidal, Simpson's, Simpson's $3/8$th, Milne's, Boole's, and other higher order rules. As we will explain later, our method makes the largest difference for the rectangular and midpoint rule because it is in these cases that $N'=N$. Thus, we will primary focus on these rules and denote them as $\mathcal{Q}_R$ and $\mathcal{Q}_M$, respectively. Because the interval is evenly partitioned, $\Omega_i=\frac{1}{N'}=\frac{1}{N}$, the integral estimate for $\mathcal{Q}_R$ is,
\begin{align}
\hat{I}_{\mathcal{Q}_R}= \frac{1}{N}\sum_{i=1}^{N}f(\vec{x}_i),
\end{align}
and for $\mathcal{Q}_M$ it is,
\begin{align}
\hat{I}_{\mathcal{Q}_M}= \frac{1}{N}\sum_{i=1}^{N}f(\frac{\vec{x}_{i+1}+\vec{x}_i}{2}).
\end{align}
The positions used to evaluate the function are known ahead of time, which is informationally different than the randomly generated samples in MC.




Rather than computing a standard deviation as a measure of fidelity, the fidelity of quadrature integration methods are represented in the literature by the error,
\begin{align}
\epsilon_{\mathcal{Q}}=I-\hat{I}_{\mathcal{Q}}.\label{error}
\end{align}
By assuming that $f$ may be approximated by a $D$ dimensional piecewise spline interpolant per containing volume $\Omega_{i}$, to leading order the error is,
\begin{align}
\epsilon_{\mathcal{Q}}\approx\sum_{i=1}^{N'}\frac{\sum_{d=1}^Df_d^{(n_{\mathcal{Q}})}(c_{\mathcal{Q},i})}{C_{\mathcal{Q}}{N'}^{\frac{n_{\mathcal{Q}}}{D}+1}}.\nonumber
\end{align}
The terms with a subscript ``$\mathcal{Q}$" are quantities that depend on the quadrature method implemented: $n_{\mathcal{Q}}$ is the number of derivatives (along a single direction of one of the dimensions),  $c_{\mathcal{Q},i}$ is the $i$th location the $n_{\mathcal{Q}}$th derivative of the function is considered, and $C_{\mathcal{Q}}$ is a constant.  
Letting $\overline{f}^{(n_{\mathcal{Q}})}=\frac{1}{N'D}\sum_{i=1}^{N'}\sum_{d=1}^Df_d^{(n_{\mathcal{Q}})}(c_{\mathcal{Q},i})$ be the average (which accounts for cancellations), is the error term of the Newton--Cotes formula,
\begin{align}
\epsilon_{\mathcal{Q}}\approx\frac{D\overline{f}^{(n_{\mathcal{Q}})}}{C_{\mathcal{Q}}{N'}^{\frac{n_{\mathcal{Q}}}{D}}}.\label{QuadError}
\end{align}
 It should be reiterated that this analysis is only valid for integrand functions that are well approximated by these integrated piecewise spline interpolants.
 
For $\mathcal{Q}_{R}$ and $\mathcal{Q}_M$,
\begin{align}
\epsilon_{\mathcal{Q}}\approx\frac{D\overline{f}^{(n_{\mathcal{Q}})}}{C_{\mathcal{Q}}{N}^{\frac{n_{\mathcal{Q}}}{D}}},\label{errorRM}
\end{align}
where for rectangular $n_{\mathcal{Q}_R}=1$ and for midpoint $n_{\mathcal{Q}_M}=2$.
The error of these quadrature methods decrease $\propto N^{-\frac{n_{\mathcal{Q}}}{D}}$ which instead does depend on the dimension of the integral $D$.

Naively these results suggest that MC should outperform quadrature when $D$ is greater than $(2,4)$, for $(\mathcal{Q}_R,\mathcal{Q}_M)$, respectively, given the dependence on $N$ in equations (\ref{stdMC}) and (\ref{errorRM}). Similar comparisons may be made between MC and the other quadrature methods by finding their respective $N'=N'(N)$ functionality and substituting into (\ref{QuadError}). These comparisons are ultimately naive because error and standard deviation are not equivalent measures of fidelity, the polynomial form information of $f$ is used in quadrature but not for MC, and due to the implications of comparing error to $\sigma_f$ under condition \ref{def}. 

\section{Quadrature Variance\label{quadfid}}

In this section the variance of simple quadrature and MC are compared directly by enforcing informational consistency using CS constraints. 

\subsection{A probabilistic equivalence through shuffling}


We motivate a probabilistic equivalence between random and deterministic sampling methods using an analogy where an ordered deck of cards
ends up playing the role of the function to be sampled:

``Consider an ordered deck of playing cards, $f$, with infinitely many cards, $|\Omega|$. Random drawing (MC sampling) from the ordered deck is drawing cards from random or unknown positions in the deck. 
This is MC$[f]$. Now consider shuffling the deck completely and randomly such that $f\stackrel{\mathcal{S}}{\rightarrow}\mathcal{S}(f)= f'$. Because the deck is completely shuffled, the positions of the cards are randomly placed. This means that it does not matter if you were to draw cards off the top of the deck, evenly spaced throughout the deck $\mathcal{Q}[f ']$ (deterministic quadrature grid sampling), or if you were to draw them from random positions in the deck MC$[f ']$ -- the draws are just as random. 
We find that it is self evident that quadrature drawing from a shuffled deck is probabilistically equivalent to randomly drawing cards from an ordered deck."
More compactly we can represent this equivalence as\footnote{Without this equivalence casinos would go out of business.}
\begin{align}
\mathcal{Q}[f']=\mbox{MC}[f]=\mbox{MC}[f'].\label{MCQ}
\end{align}


We seek to equate drawing cards at random to MC sampling; however, we need to be careful in our approach. Drawing cards is done without replacement where MC function sampling is usually done with replacement for convenience. We will first express MC in discrete terms and then take the continuous limit to be rigorous.


Because we are using MC samples to estimate an \emph{integral}, in general the draws should be made \emph{without} replacement in theory or else regions of the integral might get double counted, which would bias the result \cite{MCOHagan}. This means the joint probability distribution for drawing a collection of function values from $N$ samples of a discrete function follows a multivariate hypergeometric distribution, as we show.

The common heuristic analogy for hypergeometric distributions are that they describe the statistics of drawing marbles from an urn without replacement. We will instead use the analogy of drawing cards from a shuffled deck of cards without replacement as the statistics are equivalent to marble drawing without replacement. The multivariate hypergeometric distribution that represents sampling for discrete MC is,
\begin{align}
p_{\mbox{\tiny{MC}}}(N_1,...,N_C|K_1,...,K_C,N,K)=\frac{\prod_{c=1}^C\binom{K_c}{N_c}}{\binom{K}{N}},\label{hyperMC}
\end{align}
 where $c$ is the function value (card value) of the sample (draw), $C=|\{c\}|$ is the total number of distinct function values (distinct card values), $N_c$ is the number of samples (draws) having function value $c$, $K_c$ is the total number of function values (cards) having value $c$ (in the deck), $K=\sum_{c=1}^CK_c$ is the total number of function values (cards), and $N=\sum_{c=1}^C N_c\leq K$ is the total number of samples (draws). If instead one sampled from random positions \emph{with} replacement, then the statistics above would instead follow a multinomial distribution, which is i.i.d. 

The hypergeometric picture of a discretized MC has the correct continuous limits. The expected value of $f$ from $N$ draws,
\begin{align}
\expt{f}_{K,C}=\sum_{\{N_c\}}p_{\mbox{\tiny{MC}}}(N_1,...,N_C|K_1,...,K_C,N,K)\Big(\frac{1}{N}\sum_{c=1}^Cf_cN_c\Big)=\sum_{c=1}^Cf_c\frac{K_c}{K},\nonumber
\end{align}
 is the integral of interest in the $(K,C)\rightarrow\infty$ continuous limit because $f_c\rightarrow f$, $\frac{K_c}{K}\rightarrow \rho_{\mbox{\tiny{MC}}}(f)\,df$, which is,
\begin{align}
    I=\expt{f}=\int_{f\in \{c\}}f\cdot\rho_{\mbox{\tiny{MC}}}(f)\,df.
\end{align}
The integral is represented in the function space variables rather than the position space variables. We can switch the representation by making the objective substitution,
\begin{align}
\rho_{\mbox{\tiny{MC}}}(f)=\int_{\vec{x}\in\Omega}\rho(f|\vec{x})\rho(\vec{x})\,d\vec{x}=\frac{1}{\Omega}\int_{\vec{x}\in\Omega}\delta(f-f(\vec{x}))\,d\vec{x},\label{3.9}
\end{align}
above, and then integrating over $f$ to obtain,
\begin{align}
I=\int_{f\in \{c\}}f\cdot\Big(\frac{1}{\Omega}\int_{\vec{x}\in\Omega}\delta(f-f(\vec{x}))\,d\vec{x}\Big)\,df=\frac{1}{\Omega}\int_{\vec{x}\in\Omega}f(\vec{x})\,d\vec{x},
\end{align}
which demonstrates the equivalence of the function space and position space representations. The variance of the expected value of $f$ for $N$ draws out of $K$ is,
\begin{align}
\mbox{Var}_{K,C,MC}\equiv \mbox{Var}_{K,C,MC}\Big(\frac{1}{N}\sum_{c=1}^Cf_cN_c\Big) =\frac{K-N}{N(K-1)}\Big(\sum_{c=1}^Cf_c^2\frac{K_c}{K}-\expt{f}_{K,C}^2\Big),\nonumber\\
\label{MCvarhyper}
\end{align}
which in the $(K,C)\rightarrow\infty$ continuous limits is the anticipated result, $\sigma^2_{\hat{I}_{\mbox{\tiny{MC}}}}=\sigma_{f}^2/N$, using the same arguments as before. Because the factor $\frac{K-N}{K-1}\rightarrow 1$ in the continuous limits, it does not matter if the sampling is done with or without replacement in practice due to the continuity of $\Omega$. 

\begin{figure} 
\includegraphics[width=\textwidth]{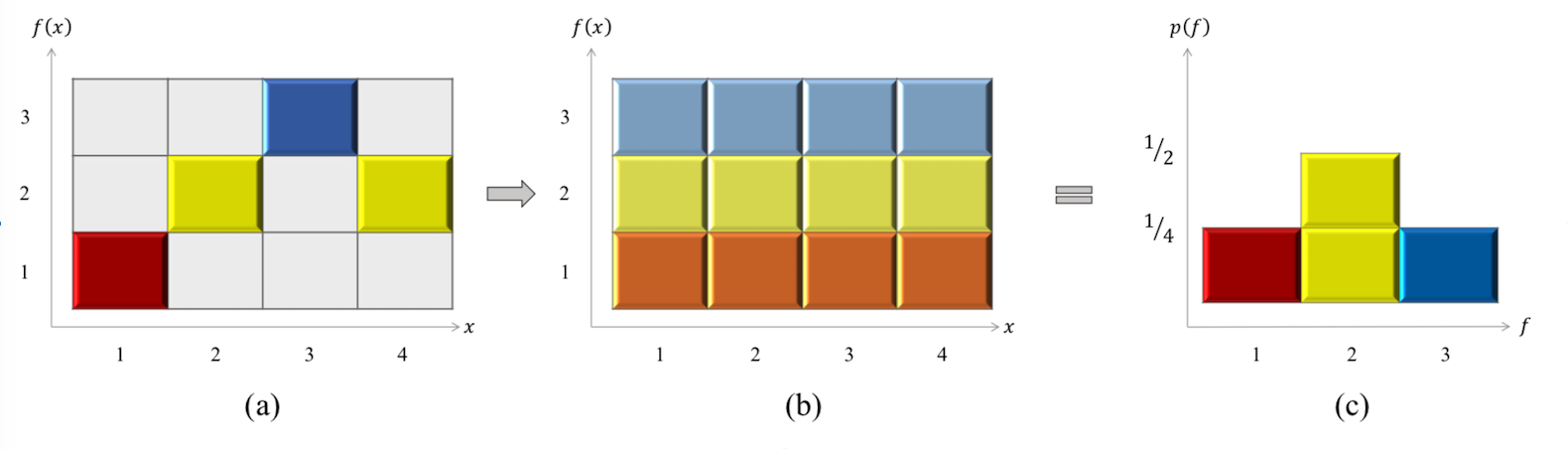}
\caption[]{This figure depicts a discrete example of MC and quadrature equivalence under function shuffling. Figure 1 (a): This is a discrete function $f(x)$ that takes the values $f$ of $(1,2,3,2)$ when evaluated at $x$ locations $(1,2,3,4)$, respectively. Drawing a single random sample from this function results in the probability depicted in Figure 1 (c). Figure 1 (b): The effect of shuffling (a) results in a probabilistically smeared function across $x$ where positional dependence is lost (hence the dimmer shades of color). Shuffling causes quadrature sampling to become position independent. The resulting probability distribution of a single quadrature sample from Figure 1 (b) is thus also given by Figure 1 (c). }
\label{Figure1}
\end{figure}

 This explicit analysis of MC has showed that drawing from a shuffled deck of cards is an identical representation of MC integrand sampling and vice versa. Thus, the probabilistic equivalence (\ref{MCQ}) holds for a deck of cards or for MC, which is what we desired to express. An example of equivalence (\ref{MCQ}) for discrete function sampling and shuffling is given in Figure \ref{Figure1}.

 \subsubsection{Shuffling Quadrature over the same integrand as MC}


Let's consider the discrete setting for quadrature integration and what it would mean to impose that it is objectively consistent with MC. The only way to guarantee that the quadrature and MC integration methods are objectively consistent is if they are guaranteed to be over the same integrand function $f$. In the discrete setting, this means that both MC and quadrature would need to at least agree on the set of values $\{K_c\}$; however, there may still be epistemic uncertainty in the way that the values of $f_c$ are distributed in $\Omega$. If it was not the case that both methods agree on $\{K_c\}$, then quadrature could assign nonzero probabilities to values that are physically impossible to produce from $f(x)$ due to a lack of objective constraints on the subjective Bayesian priors. 

A general discrete quadrature sampling probability distribution that is objectively consistent with MC is therefore,
\begin{align}
p_{\mathcal{Q}}(N_1,...,N_C|K_1,...,K_C,N,K,\vec{g},\mathcal{I}),\label{hyperquad}
\end{align}
where $\vec{g}=(\vec{x}_1,...,\vec{x}_N)$ denotes the sampling grid and $\mathcal{I}$ pertains to the prior information about the joint function value locations. Thus, rather than being a ``draw $N$ cards randomly without replacement from the deck", the drawer knows ahead of time from where in the deck they will draw cards and the probabilities of their outcomes before sampling as well as potentially some information about values not sampled.

Because nonuniform grid location sampling is isomorphic to uniform grid location sampling of an appropriately reordered function, grid location marginalization represents the effect that shuffling has on quadrature sampling probabilities, i.e. marginalizing over the sample configuration space is the same as shuffling. The result of shuffling the quadrature method must equal the MC probability from (\ref{hyperMC}) due to (\ref{MCQ}) and the agreement about $\{K_c\}$. If one has prior information $\mathcal{I}$, then marginalizing over the distribution $p(\vec{g})=\frac{1}{G}$, where $G$ is all possible joint sampling locations (uniformly distributed in space or not), must yield a result that is equal to (\ref{hyperMC}),
\begin{align}
&p_{\mathcal{Q}}(N_1,...,N_C|K_1,...,K_C,N,K,\mathcal{I})\nonumber\\
&=\sum_{\vec{g}\in G}p_{\mathcal{Q}}(N_1,...,N_C|K_1,...,K_C,N,K,\vec{g},\mathcal{I})\,p(\vec{g})\nonumber\\
&=p_{\mbox{\tiny{MC}}}(N_1,...,N_C|K_1,...,K_C,N,K),\label{hyperequal}
\end{align}
due to $(\ref{MCQ})$, where $d\vec{g}=\prod_{i=1}^Nd\vec{x}_i$.  The prior information $\mathcal{I}$ that represents one's knowledge of the \emph{positions} of values of the function to be integrated, is washed out by shuffling. This equality and analysis holds even in the more general case where one also marginalizes over partial knowledge of the function's spectrum $\{K_c\}$.

In the continuous limit, this equation is,
\begin{align}
\rho_{\mathcal{Q}}(f_1,...,f_N|\mathcal{I})
&=\int_{\vec{g}\in\Omega^N}\rho_{\mathcal{Q}}(f_1,...,f_N|\vec{g},\mathcal{I})\rho(\vec{g})\,d\vec{g}&\nonumber\\
&=\rho_{\mbox{\tiny{MC}}}(f_1)...\rho_{\mbox{\tiny{MC}}}(f_N),&\label{hyperequalc}
\end{align}
 which again imposes that the Bayesian form of $\rho_{\mathcal{Q}}(f_1,...,f_N|\vec{g},\mathcal{I})$ is objectively constrained to the function and the MC distributions under function or grid shuffling. We have suppressed the continuous value function frequencies $\{K_c\}\rightarrow \{K_f\}$ for notational convenience.   
 
\subsubsection{Counterfactual Shuffling (CS)}

To make direct and fair comparisons between quadrature sampling methods and MC sampling methods, we enforce what we call the CS constraints on the quadrature sampling procedure. The CS constraint is that, ``\emph{if} the function were shuffled or randomized, then it satisfies equation (\ref{hyperequalc})".  It should be emphasized that satisfying the CS constraints does not require shuffling or randomization of any kind as the shuffling condition is purely counterfactual. The CS constraints are satisfied for quadrature methods that are conditioned on the same objective function having spectrum $\{K_c\}$. This constrains the possible forms of Bayesian priors that could be put on the function sampling method to be only those that are informational consistent with MC for the purpose of direct comparisons.
  


\subsection{Quadrature variance when $f$ at different positions is unknown\label{section3.1}}

We impose the CS constraints in the case in which the value of the function at different positions is unknown and derive the consequences.
Because quadrature methods do not use randomization, the only time quadrature methods have equal information content prior to sampling as MC is when the functional form of $f$ is unknown, denoted $f_u$. In such a case, it does not matter if $f_u$ is sampled randomly or not because there is no positional information available to condition, i.e. $\mathcal{I}=\emptyset$ is the empty set. This is the case of drawing cards from a deck where one does not know the order of the cards in the deck, and therefore it doesn't matter if one shuffles the deck beforehand or not.

Knowing the quadrature locations ahead of time $\vec{g}=(\vec{x}_1,...,\vec{x}_N)_{\mathcal{Q}}$ does not inform one about the expected outcomes of $\{f_{u_i}\}_{\mathcal{Q}}$ because there is no known functional order and thus the joint probability in (\ref{hyperequalc}) is independent of $\vec{g}$,
\begin{align}
\rho_{\mathcal{Q}}(f_1,...,f_N|\vec{g},\mathcal{I}=\emptyset)\rho(\vec{g})=\rho_{\mathcal{Q}}(f_1,...,f_N)\rho(\vec{g})
\end{align}
Thus, the independence of $\vec{g}$ forces the probability distributions to be equal for simple quadrature and MC in this case,
\begin{align}
\rho_{\mathcal{Q}}(f_1,...,f_N|\vec{g},\mathcal{I}=\emptyset)=\rho_{\mathcal{Q}}(f_1,...,f_N)=\rho_{\mbox{\tiny{MC}}}(f_1)...\rho_{\mbox{\tiny{MC}}}(f_N),\label{hyperfu}
\end{align}
due to the CS constraint and using (\ref{hyperequalc}). The relationship holds for all $\{K_f\}$ ($\{K_c\}$ in the discrete case) and therefore it also holds for any potential marginalizations over them.

The result (\ref{hyperfu}) from (\ref{hyperequalc}) also expresses that, prior to sampling, that the act of physically shuffling $f_u\rightarrow f_u'=\mathcal{S}(f_u)$ does not affect the expectation values and probabilities because they are \emph{already} position independent. These relationships may be expressed compactly as,
\begin{align}
\mathcal{Q}[f_u]=\mathcal{Q}[f_u']=\mbox{MC}[f_u]=\mbox{MC}[f],\label{fuisMC}
\end{align}
which means $\mathcal{Q}[f_u]=\mbox{MC}[f]$ whether or not $f_u$ is shuffled. 

The standard deviation of the quadrature estimate in this case can be calculated by noting that $\hat{I}_{\mathcal{Q}}$ from (\ref{linearquad}) is now a linear combination of $N$ random variables $\{f_{u_i}\}_{\mathcal{Q}}$. The standard deviation of this quantity is therefore, 
\begin{align}
\sigma_{\hat{I}_{\mathcal{Q}}}=\sqrt{\expt{(\frac{1}{N'}\sum_{i=1}^Nq_if_{u_i})^2}-\expt{\frac{1}{N'}\sum_{i=1}^Nq_if_{u_i}}^2}=\frac{\sigma_{f}}{{N'}}\sqrt{\sum_{i=1}^Nq_i^2},\label{3.1}
\end{align}
because $\sigma_{f_u}=\sigma_{f}$ from (\ref{hyperfu}) and (\ref{fuisMC}). The minimum variance of the estimate occurs when $q_i=\frac{N'}{N}$ for all $i$, which gives min$(\sigma_{\hat{I}_{\mathcal{Q}}})=\sigma_{\hat{I}_{\mbox{\tiny{MC}}}}$.

For the rectangular and midpoint quadrature rules, $\mathcal{Q}(i,\{f_{u_1},...,f_{u_N}\})\rightarrow f_i$, $N'=N$, and the weights are all equal to $q_i=1$ for the $D$ dimensional integral estimate. Therefore,
\begin{align}
\sigma_{\hat{I}_{\mathcal{Q}_R}}=\sigma_{\hat{I}_{\mathcal{Q}_M}}=\frac{\sigma_{f}}{\sqrt{N}},\label{quadisMC}
\end{align}
which is equal to MC. 
When the functional form of $f$ is unknown, the rectangular rule and the midpoint rule have an identical fidelity to MC, which itself has fidelity estimates that are impartial to the knowledge of the functional form of $f$, i.e. the right hand side of (\ref{fuisMC}). In this sense the methods are statistically the same and share condition \ref{def}. Because integration is a sum in the infinite limit, $\hat{I}_{\mathcal{Q}}$ may be considered to be an unbiased estimator.

The maximum value of $\sum_{i=1}^Nq_i^2$ occurs when a single $q_{i'}=N'$ and the rest are equal to zero. In this case, effectively, a single $f_{i'}$ is used to estimate the integral and all of the other $f_i$ values are thrown away. This indeed results in the logical upper bound $\frac{\sigma_{f}}{\sqrt{N}}\leq\sigma_{\hat{I}_{\mathcal{Q}}}\leq \sigma_f$.

The values of the weights for trapezoidal, Simpson's, and the other quadrature rules are dimension dependent and unequally weighted. 
For most of these rules, and especially when $D$ is large, only a few of the evaluations may end up carrying the majority of the weight and thus many of the evaluations are effectively thrown away. Thus these methods behave worse than rectangular, midpoint, and MC when the function is not known. It should be noted that these results do not render other quadrature integration methods useless; rather, it is proven that when integrating functions with unknown functional form, there is no reason to weight some samples more than others. With the right pieces of information about the functional form of $f$, these methods are still useful.


After performing quadrature integration, one can calculate population averages and population standard deviations if desired; however, in practice, the function should be $f_u$ and thus \ref{def} persists. This evades common arguments against the use of quadrature integration that involve knowing worst case functional form information, like fine-tuned function periodicity or dependency, that would force the population averages and standard deviations to greatly deviate from the ``counterfactually known statistics", because, in fact, the functional order is \emph{not} known in this case by definition. 

\paragraph{Example:}
Here we give an example to further clarify (\ref{quadisMC}) and show how this result could have been anticipated. 

The result in question holds specifically for computational integrations over functions where the values of $f$ at different positions is \emph{unknown} as well as their positional interrelationships. To test these results without bias, one therefore cannot choose a given deterministic function with a known analytic form because such a case implies that indeed one \emph{does} have information about $f(x)$ and its potential positional correlations (perhaps through the knowledge of function continuity).\footnote{One can imagine the situation in which a friend chooses a deterministic function for me to integrate without telling me its functional form. However, because my friend knows the form of the function, the study is biased as the friend ends up choosing a function that may be arbitrarily good for quadrature integration (or bad or in-between).   } For this reason, we need considerable uncertainty in the function space to avoid function selection bias, i.e., the performance of either method in general depends entirely on the function we are testing so we remove this bias with uncertainty.

One can model the integration of a deterministic function with ``a completely unknown functional relationship" as the integration over a function that is generated from an independent and homogeneous uncertain process (with the additional stipulation that sampling the same location gives the same result). One models such a function as being generated from,
\begin{eqnarray}
f(\vec{x})\sim \rho(f|\vec{x})=\rho(f)\label{neweq}
\end{eqnarray}
which denotes that when $\vec{x}$ is sampled, that $f(\vec{x})$ is drawn from $\rho(f)$ -- a distribution independent of $\vec{x}$. After sampling at $\vec{x}=\vec{x}_i$ and observing $f_i$, this distribution collapses at $\vec{x}_i$,
\begin{eqnarray}
f(x_i)\sim\delta(f(\vec{x}_i)-f_i),
\end{eqnarray}
but remains equal to (\ref{neweq}) everywhere else. Note that we intentionally left the definite value of the integral unknown as this is commonly the case in practical estimation -- the average integral of the to be sampled function is $\mu_f=\int f\cdot\rho(f)df $ with standard deviation $\sigma_f$. Once the function has been sampled $N$ times, the only constraint on the unknown underlying function $f(\vec{x})$ is that it is known to pass through the coordinates $\{f_1,...,f_N\}$ at the respective sampled coordinates $\{\vec{x}_1,...,\vec{x}_N\}$. 

\begin{figure} 
\includegraphics[width=\textwidth]{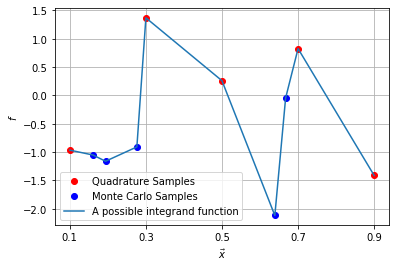}
\caption[]{We plot 5 midpoint quadrature samples in red, and 5 MC samples in blue. The light blue line connecting them is one of many possible realizations of the true underlying unknown integrand function. The true underlying unknown integrand function must pass through both the quadrature and MC samples for a given integrand function due to the CS constraints.  \label{revisedfigure} }
\label{Figure2}
\end{figure}

In Figure \ref{revisedfigure} we show the results of a simple comparison experiment after a few samples. We let $\rho(f)$ be a normal distribution with zero mean and variance equal to one. Because knowing the location of the samples ahead of time does not provide any additional information, a sample at any $\vec{x}$ is just as random in $f$ as any other sample -- this provides the intuition behind (\ref{fuisMC}) and (\ref{quadisMC}). We quadrature sample and MC sample our unknown yet deterministic function 100 times each and compare their population average and standard deviations. As expected, the integral estimations and population standard deviations closely follow the expected values (trivially) and they have no dependence on the dimension of the input space. Repeating this 10,000 times is considering the integration of 10,000 different integrand functions. Again we find average population integral statistics that closely match one another and the theory. The difference between the computation complexity of the methods is attributed to the difference in time it takes to generate $N$ random MC samples as compared to $N$ quadrature samples, which is usually small in practice. MC methods in general become slower as the quality of the random samples improve.

\subsection{Quadrature fidelity when one has position information about $f$\label{quadfidsection}}
 We seek to answer what we feel is a reasonable claim: ``How could having more information about $f$ possibly lead to situation in which the fidelity of a quadrature estimate is any worse, beyond correction or contextual qualification, than if we did not have this information"? We again impose the CS constraints.

Equations (\ref{hyperequal}) and (\ref{hyperequalc}) can be used to generate fidelity relationships between simple quadrature and MC. The variance of $f$ given $N$ draws out of $K$ is (\ref{MCvarhyper}), which is computed using $p_{\mbox{\tiny{MC}}}$ from (\ref{hyperequal}). The quadrature variance of the integral estimator conditional on: $N$, $K$, $\vec{g}$, and $\mathcal{I}$, is computed with (\ref{hyperquad}) and is a conditional variance,
\begin{align}
\mbox{Var}_{K,C,\mathcal{Q}|\vec{g}}\equiv \mbox{Var}_{K,C,\mathcal{Q}}\Big(\frac{1}{N}\sum_{c=1}^Cf_cN_c\Big|(\vec{g},\mathcal{I})\Big),
\end{align}
Due to (\ref{hyperequal}) and (\ref{hyperequalc}), we can use the law of total variance to formulate a relationship between the variances of quadrature and MC,
\begin{align}
\mbox{E}[\mbox{Var}_{K,C,\mathcal{Q}|\vec{g}}]=\sum_{\vec{g}\in G}p(\vec{g})\,\mbox{Var}_{K,C,\mathcal{Q}|\vec{g}}
\leq \mbox{Var}_{K,C,MC}.\label{result}
\end{align}
If only partial information is known about the function spectrum, $\{K_c\}$, then the unknown parts can also be marginalized over, which again leads to an analogous result. In the extreme case in which $\{K_c\}$ is completely unknown and thus it is completely marginalized over or $\mathcal{I}=\emptyset$, the inequality becomes an equality as the joint probability becomes independent of $\vec{g}$. Using analogous arguments from the previous subsection, this expression holds exactly for $\mathcal{Q}_R$ and $\mathcal{Q}_M$ because $N'=N$, that is, 
\begin{align}
\mbox{E}_{\rho(\vec{g})}[\sigma_{\hat{I}_{\mathcal{Q}_R}}]=\mbox{E}_{\rho(\vec{g})}[\sigma_{\hat{I}_{\mathcal{Q}_M}}]\leq \sigma_{\hat{I}_{\mbox{\tiny{MC}}}}=\sigma_f/\sqrt{N},\label{result2}
\end{align}
for arbitrary $f$, which motivates the title of this article.

What we have found is that \emph{if} you were to randomly shuffle the function, i.e. effectively performing MC random sampling by equation (\ref{hyperequalc}), then it is expected that the variance increases. An alternative interpretation of this inequality is that the quadrature variance of an arbitrary integrand function with arbitrary prior information $\mathcal{I}$ is on average less than or equal to its corresponding MC variance (given CS constraints). Furthermore, because $\int\rho(\vec{g})\hat{I}_{\mathcal{Q}|\vec{g}}\,d\vec{g}=I$, it means that on average the CS constrained quadrature method is an unbiased estimator of the integral. Further, the infinite sample limit is the continuous limit Riemann sum definition of integration. 


In Appendix \ref{appendix2} we make brief comments and comparisons to other MC type integration methods including Markov Chain Monte Carlo (MCMC), importance sampling, stratified sampling, and Latin Hypercube sampling; however, a direct comparison to some of these methods are outside the scope of this article.

\paragraph{Example:}

\begin{figure} 
\includegraphics[width=\textwidth]{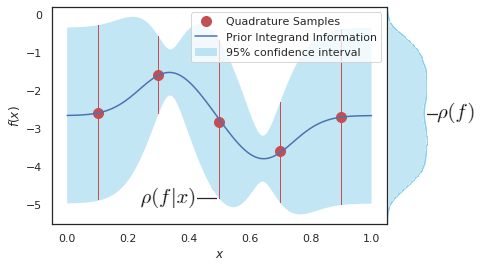}
\caption[]{We plot 5 midpoint quadrature samples in red, their $95\%$ confidence intervals with red lines, and the remaining confidence intervals of $\rho(f|\vec{x})$ in light blue. The blue line is one of many possible realizations of the true underlying unknown integrand function, analogous to Figure \ref{revisedfigure}. On the right hand side of the plot is the marginalized distribution $\rho(f)$, which is the MC distribution, and is found by marginalizing $\rho(f|\vec{x})$ uniformly over $\vec{x}$ (i.e. shuffling). Visually inspecting the length of the red lines and averaging them allows one to anticipate inequality (\ref{result2}). The notion of correcting expected errors is discussed at the beginning of Section \ref{sectionerror}.   \label{quadmix} }
\label{Figure3}
\end{figure}

Analogous to our previous example, we aim to remove function selection bias by modeling the CS constrained deterministic function with uncertainty. A deterministic integrand function with a completely known or partially known functional relationship" is represented with a function that is to be generated from heterogeneous uncertain process (again with the stipulation that sampling the same location gives the same result). One models such a function as being generated from,
\begin{eqnarray}
f(\vec{x})\sim \rho(f|\vec{x}),\label{neweq2}
\end{eqnarray}
which denotes that before $\vec{x}$ is sampled, that $f(\vec{x})$ is drawn from $\rho(f|\vec{x})$. After sampling at $\vec{x}=\vec{x}_i$ and observing $f_i$, this distribution collapses to,
\begin{eqnarray}
f(x_i)\sim\delta(f(\vec{x}_i)-f_i),
\end{eqnarray}
singularly at $\vec{x}_i$ but remains equal to (\ref{neweq2}) everywhere else. In Figure \ref{quadmix} we demonstrate how the result (\ref{result2}) can be made more intuitive.  

\subsubsection{Discussion}

By finding a way to directly compare quadrature and MC, this analysis clarifies some differences between Bayesian and frequentist methods. On the one hand, Bayesian priors can assign nonzero probability to impossible values of the function $f(\vec{x})$ simply because the set and number of possible function values are not known \emph{a priori}; however, a Bayesian can use their prior information to calculate a \emph{value} for the anticipated error and variance prior to sampling (or try many different priors). The values of the probabilities used in the calculations are assigned by the Bayesian before sampling, are epistemic in nature, and are updatable. On the other hand, it is as if frequentists (when considering sampling statistics) interpret the probability distribution somewhat similar to a physical random process state variable. Thus a frequentist does not express the \emph{value} of the probability before sampling (and thus cannot state a prediction prior to sampling) but instead manipulates the \emph{counterfactually}  assumed ``physical" state variables $I,\sigma_f$, and $\rho_{\mbox{\tiny{MC}}}$ and finds relationships between them. However, by using these true yet unknown quantities, a frequentist is never ``wrong", by definition, because the values are assumed to be the physical random process's true values \emph{a priori}. The informational shortcomings of both inference methods are symptoms that stem from the fact that, in practice and by definition, we don't know the true value(s) of the object we are trying to estimate or make inferences about. If we did know the true value(s), then every inference would be trivial.  

Due to the counterfactual nature of condition \ref{def} and the enforcement of CS to make direct comparisons, naturally this Bayesian simple quadrature analysis becomes more frequentist in nature. 
Constraining the Bayesian inference to be consistent with the frequentist analysis upon shuffling forces the spectrum $\{K_c\}$ to be assumed to be known, similar to \ref{def}, which, as far as integration is concerned, is as informative as knowing $f(x)$ or $I$ completely. Thus in practice the types of prior information that are consistent under CS are limited to function value location knowledge and permutation.\footnote{It should be noted that something like the Bayes-Hermite Gaussian Process quadrature method \cite{BQ}, which in principle can allow the function to take any value with at least infinitesimal probability, does not fall into this category in general.} Allowing $\{K_c\}$ to be uncertain mitigates this effect somewhat as it begins to push the frequentist MC method into a more Bayesian domain as expressing and using the value of $\rho(\{K_c\})$ (a prior) before sampling requires the use of Bayesian probability theory. We can summarize our result (\ref{result2}) as follows, ``Without making any more assumptions than are already present in MC's fidelity analysis, it is expected that simple quadrature methods have fidelity estimates that are just as good or better than MC".

When prior information is known about the functional form of $f$, we suggest that one should instead seek methods that provide stronger assessments of fidelity by forgoing the counterfactual use of a definite yet unknown $\sigma_f$. When possible, we instead recommend focusing one's attention on the potentially more practical scenario of determining error quantities by utilizing one's prior knowledge directly.


\section{Using Bayesian probability to improve the Newton--Cotes error quadrature error estimates\label{sectionerror}}
In this section we improve the Newton--Cotes quadrature error analysis found in the literature. First we will make some remarks about error and error estimation in practice. 

If one knows the exact value of the error, i.e. how much the estimate is biased, then one can immediately obtain an exact estimate of $I$ by creating a new estimate,
\begin{align}
\hat{I}'=\hat{I}+\epsilon_{\mathcal{Q}}=I,\label{noerror}
\end{align}
i.e., correcting it.\footnote{This is what makes supervised machine learning tasks successful -- an error is observed and the model is trained in a way that corrects it.}
This occurs because knowing (\ref{error}) with certainty accidentally implies that both $I$ and $\hat{I}$ are known too.  

By definition in \emph{estimation}, the value of the error must be uncertain in practice and therefore the error can have nontrivial expectation values. The expected value of the error behaves similar to the error itself as far as estimation goes. If $\expt{\epsilon_{\mathcal{Q}}}\neq 0$, then the estimate is expected to be biased. Thus, the expected bias can be removed simply by using the augmented estimate,
\begin{align}
\hat{I}'=\hat{I}+\expt{\epsilon_{\mathcal{Q}}}.\label{noerror2}
\end{align} 
If, for instance, the grid points are known to overestimate the integral, then using (\ref{noerror2}) reduces the error in the estimate, given $\expt{\epsilon_{\mathcal{Q}}}$. If there is no \emph{objective} reason to believe the estimate is an over or under estimate given the information $I$ from a Bayesian point of view, then by symmetry, \emph{a priori}, the estimate is expected to be unbiased $\expt{\epsilon_{\mathcal{Q}}}=0$. In this sense, most quadrature methods are unbiased or have biases that in principle may be corrected. This can be applied to any expected bias of the quadrature estimate in the previous section. 

Because known error can be corrected (\ref{noerror}), we should focus our attention on calculating the standard deviation of the unbiased error rather than just the error.\footnote{Full probability considerations of the error are favorable to the standard deviations when available in practice. Ultimately, probability distributions are more informative because they can, for example, detect bimodality and other features that simple standard deviations would hide.} Incorporating functional form information into the estimate ultimately reduces the uncertainty in the error while also providing knowledge about the expected bias that can be corrected (\ref{noerror2}). To calculate the standard deviation of the quadrature error it requires one to be more specific about what is and is not known about the functional form of $f$. We give an example here.

\subsection{Quadrature error variance due to bounded derivative information}

Let's assume that it is known that the function $f$ has a continuous $n_{\mathcal{Q}}$th derivative in $\Omega_i$ and that it is well approximated by a spline of that order. Further, it is known that max$(|f^{(n)}(\vec{x})|)\leq\eta$ in any direction -- this is only partial prior information so the derivatives are still largely uncertain. Using Bayesian probability theory, one can assign a prior probability distribution to these uncertain derivatives.  This implies that $\sigma_{f^{(n)}}\leq\sqrt{2}\eta$ is bounded, where equality is reached in the worse case of a $50/50$ probability of $f^{(n)}(\vec{x})=\pm\eta$. Due to the arguments of the error being uncertain and the limited derivative information we have, we can treat the uncertain derivatives $f^{(n_{\mathcal{Q}})}(c_{\mathcal{Q},i})$ as random variables with standard deviation $\sigma_{f^{(n)}_i}$. 

We uniform grid sample the integral such that the evaluations $\{f_i\}$ are calculated and are no longer uncertain, but the derivatives in between the samples are uncertain. The error is given by equation (\ref{QuadError}). 
Using this information we can approximate the standard deviation of the unbiased error. The standard deviation of the unbiased quadrature error in this case becomes,
\begin{align}
\sigma_{\epsilon_{\mathcal{Q}}}\approx \sqrt{\expt{\Big(\sum_{i=1}^{N'}\frac{\sum_{d=1}^Df_d^{(n_{\mathcal{Q}})}(c_{\mathcal{Q},i})}{C_{\mathcal{Q}}{N'}^{\frac{n_{\mathcal{Q}}}{D}+1}}\Big)^2}-\expt{\sum_{i=1}^{N'}\frac{\sum_{d=1}^Df_d^{(n_{\mathcal{Q}})}(c_{\mathcal{Q},i})}{C_{\mathcal{Q}}{N'}^{\frac{n_{\mathcal{Q}}}{D}+1}}}^2}\nonumber\\
=\frac{\sqrt{\sum_{i=1}^{N'}\sum_{d=1}^D\sigma^2_{f_{d,i}^{(n)}}}}{C_{\mathcal{Q}}{N'}^{(\frac{n_{\mathcal{Q}}}{D}+1)}}\nonumber.
\end{align}
Letting $\overline{\sigma}_{f^{(n_{\mathcal{Q}})}}^2=\frac{1}{N'D}\sum_{i=1}^{N'}\sum_{d=1}^D\sigma^2_{f^{(n_{\mathcal{Q}})}}$ be the average,
\begin{align}
\sigma_{\epsilon_{\mathcal{Q}}}\approx \frac{\sqrt{D}\,\overline{\sigma}_{f^{(n_{\mathcal{Q}})}}}{C_{\mathcal{Q}}{N'}^{(\frac{n_{\mathcal{Q}}}{D}+\frac{1}{2})}}\leq \frac{\,\sqrt{2D}\eta}{C_{\mathcal{Q}}{N'}^{(\frac{n_{\mathcal{Q}}}{D}+\frac{1}{2})}}.
\end{align}
For $\mathcal{Q}_R$ and $\mathcal{Q}_M$, $N'=N$, so the $N$ dependence is,
\begin{align}
\sigma_{\epsilon_{\mathcal{Q}}}\propto{N}^{-(\frac{n_{\mathcal{Q}}}{D}+\frac{1}{2})},\label{validateme}
\end{align}
which approaches $N^{-\frac{1}{2}}$ as the dimension of the space $D\rightarrow\infty$. 
As an illustration, for midpoint quadrature, the calculation assumes that the functional behavior about each point is a hyperbolic paraboloid with an unknown concavity in each direction that is bounded by $|f^{(2)}_d(c_{\mathcal{Q}_M,i})|\leq\eta$, which contributes to the error in the sum. Unlike MC, every quantity is known in this measure of fidelity (i.e. there is no condition like \ref{def}), but then again the analysis may not obey CS type constraints so these analyses cannot be directly compared with MC; however, they still may be loosely compared. 


The theoretical exponents of $N$ in the standard deviation of the error are simulated in Appendix \ref{appendix} for $\mathcal{Q}_R$ and $\mathcal{Q}_M$ in $1,2,4,8,$ and $16$ dimensions under their respective piecewise spline interpolant assumptions.  The simulated exponents closely match the theoretical exponents given above. 

Although the ``worst case quadrature error" decreases $\propto N^{-\frac{n_{\mathcal{Q}}}{D}}$, it also has the lowest probability of occurring due to the high likelihood of over and under estimates averaging out of the unbiased error. The worst case error occurs if either all of the $N'$ terms are an over or under estimate, which is analogous to the chances of flipping a coin $N'$ times and having the outcomes be all heads or all tails, i.e. the probability $2^{-N'+1}$, which is exponentially small. 

Having additional information about the functional form of $f$ can be used by quadrature methods to improve the inference of $I$. Having enough information about $f$ allows for the development of quick to converge adaptive sampling methods that exist in the literature, like Bayesian Quadrature \cite{BQ,Briol}.


\subsection{Quadrature error variance of higher order quadrature methods} The number of random derivatives $f_d^{(n_{\mathcal{Q}})}$ that contribute to the statistical reduction of the estimate's error under piecewise interpolant assumptions is equal to the number of partitions $N'$ rather than the number of samples $N$. Because $N'< N$ for methods other than rectangular and midpoint, the effective statistical reduction provided by canceling errors gets suppressed when $D$ is large. For instance, the relationship between the number of partitions and the number of samples for the trapezoidal rule and Simpson's rule is $N'=(N^{\frac{1}{D}}-1)^D$ and $N'=(N^{\frac{1}{D}}-2)^D$, respectively, which have a heavy dependency on $D$.\footnote{Using the trapezoidal rule, $N'=1$ for $N=2^{20}\sim 10^6$ when $D=20$. Using Simpson's rule, $N'=1$ for  $N=3^{20}\sim3.5\times 10^9$ when $D=20$.} Given we are interested in the performance of quadrature variance in the large $D$ case, the computation of quadrature variance does not provide much additional benefit for these methods in terms of $N'$. These methods will have to rely solely on the knowledge of $\sigma_{f^{(n_{\mathcal{Q}})}}/C_{\mathcal{Q}}$ being small.

\section{Conclusion}

We give two justifications in this article that provide evidence in favor of the claim that simple quadrature integration methods are just as good as MC integration methods. The first involves imposing CS constraints on the types of Bayesian priors one could use for quadrature integration inference. These constraints impose objectivity into the Bayesian analysis and allow us to compare quadrature and MC concretely by forming a statistical bridge between deterministic and random sampling. Without making any more assumptions than are already present in MC's fidelity analysis, we show that it is expected that rectangular and midpoint quadrature methods have fidelity estimates that are just as good or better than MC, that is, $\mbox{E}[\sigma_{\hat{I}_{\mathcal{Q}}}]\leq \sigma_f/\sqrt{N}$. We expect CS to be a useful tool for making other types of direct comparisons between Bayesian and frequentist methods.


This article promotes the use of the information available to the practitioner when performing numerical integration. We recognize that in practice the error of any estimate is uncertain.  This allows us to write the sum of composite errors in the Newton--Cotes error formula as sum of unbiased random variables, which leads to the second justification of our main claim. We show that the standard deviation of the unbiased quadrature error for composite rectangular and midpoint rules improves as $\propto N^{-(\frac{n_{\mathcal{Q}}}{D}+\frac{1}{2})}$, where the extra factor $N^{-\frac{1}{2}}$ persists even when the number of dimensions is large.

Our results contradict typical preconceptions given throughout the literature that MC is expected to outperform simple quadrature in high dimensions as both of the justifications we presented gain dimension independent $N^{-\frac{1}{2}}$ factors. This opens a space for practitioners to formulate new uncertainty quantification techniques based on simple quadrature methods with rigorous fidelity measures in areas that are typically dominated by MC type methods.

\appendix

\section{Appendix: Direct comparisons to other MC based integration methods?\label{appendix2}}

We attempt to make some direct comparisons between simple quadrature and Markov Chain Monte Carlo (MCMC), importance sampling, stratified sampling, and Latin Hypercube sampling using this analysis. We use comparisons of these methods to MC from literature and then compare them to simple quadrature using (\ref{result2}).

The variance of a MCMC estimate is,
 \begin{align}
\mbox{Var}_{MCMC}(\frac{1}{N}\sum_{i=1}^Nf_i)=\frac{1}{N}\sigma^2_{f}+\frac{2}{N^2}\sum_{j>i,i=1}^N\mbox{Cov}(f_i,f_j)\equiv\frac{\tau}{N}\sigma^2_{f}
\end{align}
where $\tau$ is the integrated autocorrelation time. When the samples are all drawn i.i.d., one finds $\tau=1$, which is MC. Due to the Markov chain nature of MCMC, the samples tend to have positive autocorrelation that exponentially decay $\mbox{Cov}(f_i,f_j)\sim \exp(-|j-i|/\tau)$, which means that $\tau$ tends to be greater than one, i.e. it tends that $\mbox{Var}_{MC}\leq \mbox{Var}_{MCMC}$. Using (\ref{result2}) this implies that the expected variance of rectangular and midpoint quadrature is also less than the variance of MCMC.

The other sampling methods are difficult to compare directly with simple quadrature for a number of reasons. The variance of importance sampling in the general case is not directly compared to MC in the literature because the ``importance sampling distribution" $h(\vec{x})$ is chosen somewhat arbitrarily in practice. It is known that the optimal choice of $h(\vec{x})$ is $h^*(\vec{x})=|f(\vec{x})|/\expt{|f(\vec{x})|}$; however, this is self defeating as it requires effectively knowing the value of the integral in the first place, $\expt{|f(\vec{x})|}$. This causes practitioners to rely on heuristics or guess and check strategies for choosing $h(\vec{x})$. Thus, by not having a rigorous general comparison to MC, it also lacks a rigorous comparison to simple quadrature.  Latin Hypercube sampling (LHS) and stratified sampling (SS) have been compared with MC, but because $\mbox{Var}_{LHS}\leq\mbox{Var}_{MC}$ (for monotonic functions at least) and $\mbox{Var}_{SS}\leq\mbox{Var}_{MC}$ while $\mbox{E}[\mbox{Var}_{\mathcal{Q}}]\leq\mbox{Var}_{MC}$, one cannot determine from these inequalities alone whether or not simple quadrature is less than the other methods or vice-versa -- for general $f(\vec{x})$. Making general direct comparisons with these methods is outside the scope of the present article.

\section{Appendix: Simulation Results}\label{appendix}

For $\mathcal{Q}_R$ and $\mathcal{Q}_M$ we simulate the standard deviation of the error under their respective piecewise spline interpolant assumptions. Using $N=2^{16}=65536$ samples, we simulated the approximate error of the integral using,
\begin{align}
\hat{\epsilon}_{\mathcal{Q}}\approx \frac{1}{C_{\mathcal{Q}}N^{\frac{n_{\mathcal{Q}}}{D}+1}}\sum_{i=1}^N\Big(\sum_{d=1}^Df_d^{(n_{\mathcal{Q}})}(c_{\mathcal{Q},i})\Big),
\end{align}
where $\{f_d^{(n_{\mathcal{Q}})}(c_{\mathcal{Q},i})\}$ are randomly generated from a uniform distribution having $\sigma_{f^{(n_{\mathcal{Q}})}_d}$ (these distributions should be assigned based on the prior knowledge and the choice of its functional form does not affect our result much -- here uniform was chosen \emph{a priori}). The error is randomly simulated 100 times $\{\hat{\epsilon}_{\mathcal{Q}}^{(1)},...,\hat{\epsilon}_{\mathcal{Q}}^{(100)}\}$ and the unbiased standard deviation of the error is estimated using,
\begin{align}
\hat{\sigma}_{\hat{\epsilon}_{\mathcal{Q}}}=\sqrt{\frac{1}{100}\sum_{j=1}^{100}(\hat{\epsilon}_{\mathcal{Q}}^{(j)})^2},
\end{align}
as the expected value of the error is zero. If the expected error is not zero because $\expt{\sum_{i=1}^N\sum_{d=1}^Df_d^{(n_{\mathcal{Q}})}(c_{\mathcal{Q},i})}\neq0$, then one can use the knowledge of this bias to correct estimate (\ref{noerror2}). 

The exponent, $\chi$, in $N^{-\chi}$ from the standard deviation of the error (\ref{validateme}) is what we seek to validate in the simulation. We estimate $\chi$ from the simulation results by using,
\begin{align}
\hat{\chi}=\log_{N}\Bigg(\frac{\sqrt{D}\,\overline{\sigma}_{f^{(n_{\mathcal{Q}})}}}{C_{\mathcal{Q}}\hat{\sigma}_{\hat{\epsilon}_{\mathcal{Q}}}}\Bigg),
\end{align}
for $\mathcal{Q}_R$ and $\mathcal{Q}_M$, which in theory is $\chi=\frac{n_{\mathcal{Q}}}{D}+\frac{1}{2}$. The $\pm$ values on the simulated $\hat{\chi}$ values in Table \ref{table2} and Table \ref{table3} were generated by repeating the whole experiment 10 times and reporting their average $\hat{\chi}\pm\hat{\sigma}_{\hat{\chi}}$ of the estimates. Once one realizes that the function derivatives can be treated as random numbers, the simulation effectively becomes the exercise of demonstrating the well known result that the standard deviation of the average of a set of random numbers decreases $\propto N^{-\frac{1}{2}}$.
\begin{table}[!htbp]
\centering
\caption{$\mathcal{Q}_R$: Theory vs. Simulation of $\chi=\frac{1}{D}+\frac{1}{2}$\ }
\begin{tabular}{cccc}
\hline
Dimension $D$ & Error exponent & Theory exponent $\chi$  & Simulated exponent $\hat{\chi}$ \\\hline
1& 1.0 & 1.5& 1.500\,$\pm$\,0.007\\

2& 0.5 & 1.0& 1.004\,$\pm$\,0.006\\

4& 0.25 & 0.75 & 0.752\,$\pm$\,0.005\\

8& 0.125 & 0.625& 0.625\,$\pm$\,0.004\\

16& 0.0625 & 0.5625& 0.564\,$\pm$\,0.007\\


\hline
\label{table2}

\end{tabular}
\end{table}
\newpage








\begin{table}[!htbp]
\centering
\caption{$\mathcal{Q}_M$: Theory vs. Simulation of $\chi= \frac{2}{D}+\frac{1}{2}$\ }
\begin{tabular}{cccc}
\hline
Dimension $D$ & Error exponent & Theory exponent $\chi$  & Simulated exponent $\hat{\chi}$ \\\hline
1& 2.0 & 2.5& 2.501\,$\pm$\,0.006\\

2& 1.0 & 1.5& 1.498\,$\pm$\,0.005\\

4& 0.5 & 1.0 & 1.003\,$\pm$\,0.006  \\

8& 0.25 & 0.75& 0.748\,$\pm$\,0.006  \\

16& 0.125 & 0.625& 0.625\,$\pm$\,0.004 \\


\hline
\label{table3}
\end{tabular}
\end{table}
















\end{document}